\documentclass [10pt,a4paper,draft] {article}
\usepackage [applemac] {inputenc}
\usepackage[french] {babel}
\usepackage {amsfonts}
\usepackage {amsmath}
\usepackage {amssymb}

\begin{document}

\begin{center}
\begin{LARGE}
Notes sur la notion d'invariant caract\'eristique
\end{LARGE}
\bigskip

par : Mustapha RA\"IS  \`a Poitiers
\end{center}

\bigskip
\vspace{10 mm}

\begin{large}
\noindent
\textbf{	Introduction}
\end{large}

\medskip
	Soient $V$ un espace vectoriel de dimension finie sur un corps $\Bbb{K}$ de caract\'eristique z\'ero et
$G$ un sous-groupe du groupe lin\'eaire $GL(V)$ de $V$. Le groupe $G$ admet divers types
d'\textsl{invariants} : orbites, fonctions polyn\^omes, fonctions diff\'erentiables et distributions
(ceci lorsque $\Bbb{K} = \Bbb{R}$). En sens inverse, ayant privil\'egi\'e un ensemble d'invariants
particuliers, on peut se demander quel est le plus grand sous-groupe $\tilde{G}$ de $GL(V)$ qui
laisse invariant ces invariants particuliers. Lorsque $G = \tilde{G}$, ou au moins lorsque $G$ et
$\tilde{G}$ ont la m\^eme composante neutre (cas $\Bbb{K} = \Bbb{R}\ldots$), on pourra dire que les
invariants consid\'er\'es sont des \textsl{invariants caract\'eristiques} de $G$.

	Un autre aspect de ce probl\`eme est le suivant : Supposons que $\Bbb{K}$ soit le corps des nombres
r\'eels et que $G$ soit un sous-groupe ferm\'e (connexe) de $GL(V)$. A chaque \'el\'ement $X$ de l'alg\`ebre
de Lie $\mathfrak{g}$ de $G$ il est associ\'e un champ de vecteurs $L_X$ sur $V$ : $L_X(v) = Xv$
pour tout $v$ dans $V$, qui annule toutes les fonctions diff\'erentiables $G$-invariantes sur $V$,
lorsqu'on le consid\`ere comme un op\'erateur diff\'erentiel lin\'eaire d'ordre 1 : $(L_Xf)(v) = \, <df(v),
X.v>$. Sous l'hypoth\`ese de la connexit\'e de $G$, il vient : une fonction diff\'erentiable $f$ sur $V$ est
$G$-invariante si et seulement si : $L_Xf = 0$ pour tout $X$ dans $\mathfrak{g}$. Les invariants de
$G$ sont donc ceux de l'alg\`ebre de Lie de champs de vecteurs constitu\'ee par les $L_X$, $X$ d\'ecrivant
$\mathfrak{g}$, laquelle alg\`ebre de Lie est isomorphe \`a $\mathfrak{g}$ (on ne la distinguera pas de
$\mathfrak{g}$ \'eventuellement). On peut alors se demander, par exemple, quelle est la plus grande
alg\`ebre de Lie de champs de vecteurs \textsl{lin\'eaires} (i.e. polynomiaux, homog\`enes de degr\'e 1)
sur $V$ qui annule tous les polyn\^omes $G$-invariants, et la comparer \`a $\mathfrak{g}$ ; lorsque
cette alg\`ebre de Lie co\"incide avec $\mathfrak{g}$, on dira comme plus haut que les polyn\^omes
invariants sont caract\'eristiques pour $\mathfrak{g}$.

	Un peu plus g\'en\'eralement, notons $\mathfrak{h}$ l'ensemble des champs de vecteurs sur $V$, qui
sont combinaisons lin\'eaires des $L_X$ ($X$ d\'ecrivant $\mathfrak{g}$), \`a coefficients polyn\^omes (ou
fonctions $C^\infty$ \'eventuellement) sur $V$ ; il est clair que $\mathfrak{h}$ est une sous-alg\`ebre
de Lie de l'alg\`ebre de Lie $\chi(V)$ de tous les champs de vecteurs polynomiaux sur $V$, et que ses
invariants polynomiaux sont ceux de $\mathfrak{g}$. Comme ci-dessus, on peut chercher la
sous-alg\`ebre de Lie de $\chi(V)$ qui annule tous les polyn\^omes invariants, et la comparer \`a
$\mathfrak{h}$.

	On traite de ces questions dans la suite sous les rubriques :

\vskip 7mm
\begin{large}
\noindent
\textbf{1~Les polyn\^omes caract\'eristiques}
\end{large}

	Divers exemples. Groupes adjoints des alg\`ebres de Lie semi-simples. Alg\`ebres de Takiff. Espaces
riemanniens sym\'etriques. Alg\`ebres de Jordan.

\vskip 7mm
\begin{large}
\noindent
\textbf{2~Les orbites caract\'eristiques}
\end{large}

	Orbites caract\'eristiques dans les espaces sym\'etriques. Orbites adjointes caract\'eristiques dans les
alg\`ebres de Lie semi-simples.

\vskip 7mm
\begin{large}
\noindent
\textbf{3~Distributions caract\'eristiques}
\end{large}

	Alg\`ebres de Lie nilpotentes. Groupes unipotents.

\vfill\eject
\begin{large}
\noindent
\textbf{Quelques commentaires :} 
\end{large}

\medskip
	Sur la partie 1, il est difficile de faire preuve d'une totale originalit\'e, tant il y a de travaux sur ce type
de questions, principalement il est vrai dans le cas des groupes r\'eductifs. Le paragraphe consacr\'e aux
alg\`ebres de Takiff semble original.

	Pour l'essentiel et \`a ma connaissance, les parties 2 et 3 sont originales.

\vskip 10mm
\section{Sur les polyn\^omes caract\'eristiques}

\begin{large}
\noindent
\textbf{1.1. Divers exemples. }
\end{large}

\bigskip
	Soit $G$ un sous-groupe de $GL(V)$. On note $\Bbb{K}[V]$ l'alg\`ebre des fonctions polynomiales sur
$V$, et $\Bbb{K}[V]^G$ celle des fonctions polynomiales $G$-invariantes. Soit $G'$ le plus grand
sous-groupe de $GL(V)$ qui laisse invariants tous les polyn\^omes $G$-invariants. On a \'evidemment
$G \subset G'$.

\vskip 7mm
\noindent
\textbf{1.1.1.}~Supposons $G$ \textsl{fini}. Alors $G'$ a les m\^emes orbites que $G$, car
$\Bbb{K}[V]^G$ s\'epare les $G$-orbites ; donc $G'$ est un groupe fini (un sous-groupe de $GL(V)$,
dont les orbites sont finies, est un ensemble fini). Par ailleurs, il existe un point $v$ de $V$ (et m\^eme
une infinit\'e de tels points $v$) dont la $G$-orbite (resp. la $G'$-orbite) a m\^eme cardinal que $G$
(resp. $G'$). Donc $G=G'$.

\vskip 7mm
\noindent
\textbf{1.1.2.}~Ici $\Bbb{K} = \Bbb{R}$ et $G$ est un sous-groupe \textsl{compact} de $GL(V)$.
Pour les m\^emes raisons que ci-dessus, les $G'$-orbites sont les $G$-orbites. Par ailleurs, $G'$ est un
sous-groupe ferm\'e de $GL(V)$, qui laisse invariante une forme quadratique d\'efinie positive ; donc
$G'$ est conjugu\'e d'un sous-groupe ferm\'e du groupe orthogonal de ladite forme quadratique. Ainsi :
$G'$ est compact. 
\begin{enumerate}
\item Ici $\Bbb{K}=\Bbb{R}$, $V = \Bbb{R}^n$ et $G = SO(n,\Bbb{R})$. Alors $G' = O(n,\Bbb{R})$.

\medskip
\item Ici $\Bbb{K}=\Bbb{R},\ V=\Bbb{C}^n$ consid\'er\'e comme espace vectoriel r\'eel de dimension
$2n$, et $G = U(n)$ est le groupe unitaire. Alors $G' = O(2n)$.
\end{enumerate}

\vskip 7mm
\noindent
\textbf{1.1.3.}~Soit $H$ le groupe de Heisenberg de dimension 3 ; son alg\`ebre de Lie $\mathfrak{h}$
est engendr\'ee, comme espace vectoriel par une base $(P,Q,Z)$ avec $[P,Q]= Z$ ; on note $x,y,z$ les
coordonn\'ees sur le dual $\mathfrak{h}^*$ de $\mathfrak{h}$, relatives respectivement \`a  $P^*, Q^*,
Z ^*$, sachant que $(P^*, Q^*, Z^*)$ est la base duale de $(P,Q,Z)$. Ici le groupe $G$ qui nous int\'eresse
est l'image $Ad^*(H)$ de $H$ par sa repr\'esentation coadjointe, c'est l'ensemble des matrices de la
forme :
$$
\begin{pmatrix}
1	&0		&b\\
0	&1	&-a\\
0	&0	&1\\
\end{pmatrix} \quad (a,b) \in \Bbb{R}^2
$$

Les orbites de $G$ sont : les plans d'\'equation : $z=z_0$, avec $z_0 \not= 0$, et les points du plan
d'\'equation $z=0$. Les polyn\^omes $G$-invariants sont les polyn\^omes en $z$ :
$$
\Bbb{R}[\mathfrak{h}^*]^G = \Bbb{R}[z].
$$

Il est clair que $G'$ est le groupe des matrices :
$$
\begin{pmatrix}
\alpha	&\beta		&u\\
\gamma	&\delta	&v\\
0	&0	&1\\
\end{pmatrix}
$$
Ici les polyn\^omes $G$-invariants ne caract\'erisent pas le groupe $G$. On notera toutefois que $G$ est
caract\'eris\'e par ses orbites.

\vskip 7mm
\noindent
\textbf{1.1.4.}~Soit $G$ un groupe de Lie semi-simple complexe connexe, et soit $\mathfrak{g}$ son
alg\`ebre de Lie. Le groupe $G$ op\`ere dans $\mathfrak{g}$ par la repr\'esentation adjointe, de m\^eme
que les sous-groupes de $G$, en particulier un sous-groupe de Borel $B$ et un radical unipotent $N$.
\begin{enumerate}
\item Il est bien connu qu'on a : $\Bbb{C}[\mathfrak{g}]^B = \Bbb{C}[\mathfrak{g}]^G$ (les
invariants polyn\^omes  pour $B$ sont des invariants pour $G$) (on a un r\'esultat analogue pour les
repr\'esentations holomorphes, de dimension finie, de $G$, au lieu et place de la repr\'esentation
adjointe). On a donc : $G \subset B' = G'$, et comme $G$ est la composante neutre de $G'$ (voir
1.2.1 plus loin), il vient que $G$ est la composante neutre de $B'$.

\medskip
\item Suite \`a une v\'erification, j'avais conjectur\'e que le groupe unipotent $N$ pouvait \^etre
ca\-ract\'eris\'e par ses invariants polyn\^omes dans $\Bbb{C}[\mathfrak{g}]$ (dans le cas $G = SL(2)$, la
v\'erification est visible dans la remarque de 3.2.2.). Cette conjecture a \'et\'e d\'emontr\'ee par A. Bouaziz et
R. Yu sous la forme suivante : $G \cap N' \subset N$ ; autrement dit : la composante neutre de $N'$ est
$N$.
\end{enumerate}

\vskip 10mm
\begin{large}
\noindent
\textbf{1.2. Groupes adjoints des alg\`ebres de Lie semi-simples}
\end{large}

\bigskip 
\noindent
\textbf{1.2.1.}~Soit $\mathfrak{g}$ une alg\`ebre de Lie semi-simple complexe de dimension finie. Le
groupe $G$ auquel on s'int\'eresse ici est le groupe $\hbox{Int}(\mathfrak{g})$, i.e. le sous-groupe
connexe de
$GL(\mathfrak{g})$, d'alg\`ebre de Lie $ad(\mathfrak{g})$. Pour chaque $X$ dans $\mathfrak{g}$, on
note $L_X$ le champ de vecteurs ``adjoint'', d\'efini par : $L_X(Y) = [X,Y]$ pour tout $Y$ dans
$\mathfrak{g}$, et on d\'esigne par $\mathfrak{h}$ l'alg\`ebre de Lie des champs de vecteurs sur
$\mathfrak{g}$, qui sont combinaisons lin\'eaires des $L_X$, \`a coefficients fonctions polyn\^omes sur
$\mathfrak{g}$. On a alors le th\'eor\`eme de Dixmier (\cite{Di-1}, th\'eor\`eme 2.1).

\bigskip
\noindent
\textbf{Th\'eor\`eme} : $\mathfrak{h}$ \textit{est l'alg\`ebre de Lie des champs de vecteurs sur
$\mathfrak{g}$, qui annule tous les polyn\^omes $\hbox{Int}(\mathfrak{g}$)-invariants sur
$\mathfrak{g}$},

\smallskip
avec son :

\medskip
\noindent
\textbf{Corollaire} : (\cite{Di-1}, remarque 2.2) : \textit{L'alg\`ebre de Lie du groupe $G'$ est celle
$ad(\mathfrak{g})$ des d\'erivations int\'erieures de $\mathfrak{g}$}.

\medskip
	Autrement dit: $G$ est la composante neutre de $G'$, et conform\'ement \`a l'introduction, $G$ est
caract\'eris\'e par ses polyn\^omes invariants.

\vskip 7mm
\noindent
\textbf{1.2.2.}~Soient $G = GL(n,\Bbb{C})$ et $\mathfrak{g}= \mathfrak{g}l(n,\Bbb{C})$. Dans
\cite{Ra} (chapitre V), il est prouv\'e que le sous-groupe de $GL(\mathfrak{g})$ qui laisse invariants les
coefficients du polyn\^ome caract\'eristique est engendr\'e par $Ad(G)$ (l'image de $G$ dans
$GL(\mathfrak{g})$ par la repr\'esentation adjointe) et par la transposition des matrices. Le m\^eme
r\'esultat vaut pour $sl(n,\Bbb{C})$ (corollaire imm\'ediat du r\'esultat pr\'ec\'edent). Donc, pour
$sl(n,\Bbb{C}) = \mathfrak{g}$, le groupe $\hbox{Int}(\mathfrak{g})$ est un sous-groupe d'indice 2
du groupe $(\hbox{Int}(\mathfrak{g}))'$.

\bigskip
\noindent
\textsc{Question}: Peut-on d\'eterminer de façon pr\'ecise le groupe $\hbox{Int}(\mathfrak{g})'$,
lorsque
$\mathfrak{g}$ est une alg\`ebre de Lie simple, comme dans le cas $\mathfrak{g} = sl(n,\Bbb{C})$ ?

\vfill\eject
\begin{large}
\noindent
\textbf{1.3. Alg\`ebres de Takiff associ\'ees \`a une alg\`ebre semi-simple}
\end{large}

\bigskip 
\noindent
\textbf{1.3.1.}~On examine ici le cas des alg\`ebres de Lie $\mathfrak{g}_m$ \'etudi\'ees dans \cite{R-T}.
Soit $\mathfrak{g} = \mathfrak{g}_0$ une alg\`ebre de Lie semi-simple complexe. On d\'efinit
$\mathfrak{g}_m$ comme \'etant le quotient de l'alg\`ebre de Lie $\mathfrak{g} \otimes \Bbb{C}[T]$ par
l'id\'eal $\mathcal{J} = \sum_{k>m} \mathfrak{g} \otimes T^k$. L'alg\`ebre des polyn\^omes
$\hbox{Int}(\mathfrak{g}_m)$-invariants sur $\mathfrak{g}_m$ a \'et\'e d\'etermin\'ee dans \cite{R-T}
(th\'eor\`eme 4.5).

\vskip 7mm
\noindent
\textbf{1.3.2.}~On aura besoin d'utiliser le th\'eor\`eme de Dixmier (\cite{Di-1}, th\'eor\`eme 2.1), \'enonc\'e
pour les champs de vecteurs \textsl{\`a param\`etres} sur $\mathfrak{g} = \mathfrak{g}_0$, sous la
forme suivante :

\smallskip
	Soit $W$ un espace vectoriel complexe, de dimension finie, et soit $(e_i)_i$ une base de
$\mathfrak{g}$. On note $\partial_i$ la d\'erivation le long du vecteur $e_i$. Soit $L = \sum_i\,
\varphi_i \partial_i$, où les $\varphi_i$ sont des fonctions polyn\^omes sur $W \times \mathfrak{g}$.
Si un tel champ $L$ annule toutes les fonctions polyn\^omes $\hbox{Int}(\mathfrak{g})$-invariantes
sur
$\mathfrak{g}$, alors $L$ est une combinaison lin\'eaire des champs adjoints $L_x\ (x\in
\mathfrak{g})$, \`a coefficients dans $\Bbb{C}[W \times \mathfrak{g}]$.

\vskip 7mm
\noindent
\textbf{1.3.3.}~\textbf{Proposition} : \textit{Soit $L : W \times \mathfrak{g}_m \longrightarrow
\mathfrak{g}_m$ un champ de vecteurs qui annule les polyn\^omes
$\hbox{Int}(\mathfrak{g}_m)$-invariants. Alors $L$ est une combinaison lin\'eaire \`a coefficients dans
$\Bbb{C}[W \times \mathfrak{g}_m]$ des champs adjoints sur $\mathfrak{g}_m$}.

\bigskip
\noindent
\textsc{D\'emonstration} : Par r\'ecurrence sur $m$. Le cas $m=0$ est le th\'eor\`eme de Dixmier.

	\textbf{a)} A titre de pr\'eparation, on traite le cas $m=1$. Chaque polyn\^ome invariant $p :
\mathfrak{g} \longrightarrow \Bbb{C}$ donne 2 polyn\^omes invariants sur $\mathfrak{g}_1 :
P_1(x_0 + x_1T)= P_1(x_0,x_1)= p(x_0), \ P_2(x_0,x_1)= \, <~dp(x_0),x_1>$.

\medskip
	Soit $L = a_0 + a_1T $ un champ de vecteurs \`a param\`etres sur $\mathfrak{g}_1$, qui annule tous
les polyn\^omes $\hbox{Int}(\mathfrak{g}_1)$-invariants sur $\mathfrak{g}_1$. Il annule tous les
polyn\^omes de type $P_1$ :
$$
	0 = \, <dP_1, L>\, = \, \big ({d\over dt}\big )_0\ P_1(x+tL(x)) = \big ({d\over dt}\big )_0\, p(x_0 +
ta_0(w,x))
$$
D'apr\`es le th\'eor\`eme de Dixmier, il existe une fonction polynomiale $b_0 : W \times \mathfrak{g}_m
\longrightarrow \mathfrak{g}$ telle que :
$$
	a_0(w,x) = [b_0(w,x), x_0]
$$
\noindent
pour tous $(w,x)$ dans $W \times \mathfrak{g}_1$ (sachant que $x_0$ est tel que $x = (x_0, x_1)$).

	Ensuite, $L$ annule tous les polyn\^omes de type $P_2$ : 
\begin{eqnarray}
0 &= &<dP_2,\, L> \, = \big ({d\over dt}\big )_0\ P_2(x+tL(x)) \nonumber\\
			&= &\big ({d\over dt}\big )_0\ <dp(x_0 + ta_0(w,x)),\ x_1 + ta_1(w,x)>\nonumber\\
			&= &\big ({d\over dt}\big )_0 <dp(x_0 + t[b_0(w,x), x_0],x_1> + <dp(x_0),\, a_1(w,x)>\nonumber
\end{eqnarray}
Mais l'invariance de $p$ implique :
$$
(*)\quad \big ({d\over dt}\big )_0 < dp(x_0 + t[b_0(w,x),x_0]) = [b_0(w,x),dp(x_0)]
$$
(Noter que $dp(x_0)$ est identifi\'e ici au gradient $\nabla p(x_0)$ grâce \`a la forme de Killing de
$\mathfrak{g}$). Soit en effet un \'el\'ement $e$ de $\mathfrak{g}$, et soit $\varphi(w,x)$ une fonction
polynomiale scalaire sur $W \times \mathfrak{g}_1$. On a alors :
\begin{eqnarray}
\big ({d\over dt}\big )_0 \nabla p(x_0 + t\varphi (w,x)[e,x_0]) &= &\varphi(w,x) \big ({d\over
dt}\big )_0\
\nabla p(x_0 + t[e,x_0]) \nonumber\\
\smallskip \nonumber\\
			&= & \varphi(w,x) \big ({d\over dt}\big )_0\ \nabla p (exp(t ad e)x_0)\nonumber\\
\smallskip \nonumber\\
			&= & \varphi(w,x) [e, \nabla p(x_0)] = [\varphi (w,x)e,\, \nabla p(x_0)] \nonumber
\end{eqnarray}
Ceci d\'emontre l'\'egalit\'e (*). On a donc : 

\bigskip
\qquad \qquad $0 = \, <dP_2,L>\,  = \, <[b_0(w,x), dp(x_0)], x_1> + <dp(x_0), a_1(w,x)>$

\bigskip
\qquad \qquad $0 = \, <dp(x_0), [x_1, b_0(w,x)] + a_1(w,x)>$

\bigskip
\noindent
et ceci pour tout polyn\^ome $p$ invariant sur $\mathfrak{g}$. Toujours d'apr\`es le th\'eor\`eme de
Dixmier, il existe une fonction polynomiale $b_1 : W \times \mathfrak{g}_1 \longrightarrow
\mathfrak{g}$ telle que : 
$$
a_1(w,x) = [b_1(w,x),x_0] + [b_0(w,x),x_1].
$$
Donc : 
\begin{eqnarray}
L(w,x) &= &[b_0(w,x],x_0] + ([b_1(w,x),x_0] + [b_0(w,x),x_1])T \nonumber\\
				&=	&[b_0(w,x) + b_1(w,x)T, x_0+x_1T]\nonumber
\end{eqnarray}
le crochet pr\'ec\'edent \'etant celui de $\mathfrak{g}_1$. C'est ce qu'il fallait d\'emontrer dans le cas
$m=1$.

\bigskip
	\textbf{b)} On traite le cas g\'en\'eral de $\mathfrak{g}_m$, avec $m\geq 1$, en supposant que la
proposition est vraie pour $\mathfrak{g}_{m-1}$. On se rapportera \`a \cite{R-T} pour les d\'etails qui
sembleraient manquer ci-dessous.

\medskip
	Chaque polyn\^ome $p$ sur $\mathfrak{g}$, $\hbox{Int}(\mathfrak{g})$-invariant, fournit $(m+1)$ 
polyn\^omes $P_0,\ldots ,P_m$ sur $\mathfrak{g}_m$,  $\hbox{Int}(\mathfrak{g}_m)$-invariants :
$$
	P_j(x_0,\ldots ,x_m) = {1\over j!} \big ({d\over dt}\big )^j p(x_0+tx_1+\cdots + t^m x_m)\mid_{t=0}
$$
Or : $P_0, P_1,\ldots , P_{m-1}$ ne d\'ependent pas de la variable $x_m$, et les polyn\^omes :
$$
	Q_j(x_0,\ldots , x_{m-1}) = {1\over j!} \big ({d\over dt}\big )^j \ p(x_0 + tx_1 +\cdots + t^{m-1}
x_{m-1})\mid_{t=0}
$$
sont $\hbox{Int}(\mathfrak{g}_{m-1})$-invariants. 

\smallskip
	Soit $L = \sum^m_{k=0}\, a_k\, T^k$ un champ de vecteurs \`a param\`etres sur $\mathfrak{g}$ qui
annule les polyn\^omes $\hbox{Int}(\mathfrak{g}_m)$-invariants sur $\mathfrak{g}_m$. Lorsque $0
\leq j \leq m-1$ : 
$$
	0 = \, <dP_j, L> \, = \, <dQ_j, a_0+a_1T +\cdots + a_{m-1}T^{m-1}>
$$
Par l'hypoth\`ese de r\'ecurrence, il existe des fonctions polynomiales $b_j : W \times \mathfrak{g}_m
\longrightarrow \mathfrak{g}\ (0 \leq j \leq m-1)$ telles que : 
$$
	a_0 + a_1 T+\cdots + a_{m-1}T^{m-1} = [b_0 + b_1T+\cdots + b_{m-1} T^{m-1}, x_0+x_1 T+\cdots
+ x_{m-1}T^{m-1}]
$$
le crochet \'etant calcul\'e dans $\mathfrak{g}_{m-1}$. Notons $[\ , \ ]_m$ et $[\  , \ ]_{m-1}$ les
crochets respectifs dans $\mathfrak{g}_m$ et $\mathfrak{g}_{m-1}$, on a : 

\medskip
	\qquad $[b_0 + b_1 T +\cdots + b_{m-1} T^{m-1},\ x_0 +\cdots + x_{m-1} T^{m-1} + x_mT^m]_m$

\medskip
	\qquad $= [b_0 +\cdots + b_{m-1} T^{m-1},\ x_0 +\cdots + x_{m-1} T^{m-1}]_{m-1}$

\medskip
		\qquad $+ ([b_0,x_m] + [b_1,x_{m-1}] +\cdots + [b_{m-1}, x_1])T^m$

\bigskip
	Soit maintenant $P : \mathfrak{g}_m \longrightarrow \Bbb{C}$ une fonction polyn\^ome
$\hbox{Int}(\mathfrak{g}_m)$-invariante. On a alors :

\medskip
	\qquad $<dP, L> \, = \, <dP, [b_0 +\cdots + b_{m-1} T^{m-1}, x_0 +\cdots + x_{m-1}
T^{m-1}]_{m-1}>$

\bigskip
	\qquad $+ <dP, a_m T^m> \, = \, <dP, [b_0 +\cdots + b_{m-1} T^{m-1}, x_0 +\cdots + x_m
T^m]_m>$

\bigskip
	\qquad $- <dP, ([b_0,x_m] + [b_1,x_{m-1}] +\cdots + [b_{m-1}, x_1])T^m>$

\bigskip
	\qquad $+ <dP, a_m T^m>$,

\bigskip
\noindent
le premier terme de la derni\`ere somme \'etant nul par l'invariance de $P$. Prenons alors $P = P_m$. Il
vient, avec :

\medskip
	\hskip 30mm $P_m(x) = \, <dp(x_0),x_m> \hbox{mod} \, \Bbb{C}[x_0,\ldots , x_{m-1}]$

\bigskip
	\hskip 30mm $<dP_m, L> \, = \, <dP_m, (a_m - \sum_{0\leq j \leq m-1} [b_j, x_{m-j}])T^m>$

\bigskip
\hskip 30mm
	$<dp(x_0), a_m - \sum_j\, [b_j, x_{m-j}]> = 0$

\bigskip
\noindent
pour tout $p$, $\hbox{Int}(\mathfrak{g})$-invariant sur $\mathfrak{g}$. Il existe donc une fonction
polyn\^ome $b_m : W \times \mathfrak{g}_m \longrightarrow \mathfrak{g}$ telle que :
\begin{eqnarray}
	a_m(w,x)	&=	&\sum_{0\leq j \leq m-1}\, [b_j(w,x), x_{m-j}] + [b_m(w,x),x_0]\nonumber\\
						&=	&\sum_{0\leq j \leq m}\, [b_j(w,x), x_{m-j}]\nonumber
\end{eqnarray}
Par cons\'equent :
$$
	a_0 + a_1 T +\cdots + a_m T^m = [b_0 + b_1 T +\cdots +b_m T^m,\, x_0 +\cdots + x_m T^m]_m
$$

\vskip 7mm
\noindent
\textbf{1.3.4.}~\textbf{Corollaire} : \textbf{1)} \textit{Soit $L$ un champ de vecteurs
\textsl{lin\'eaire} sur $\mathfrak{g}_m$, qui annule les fonctions polyn\^omes
$\hbox{Int}(\mathfrak{g}_m)$-invariants sur $\mathfrak{g}_m$. Alors $L :
\mathfrak{g}_m \longrightarrow \mathfrak{g}_m$ est une \textsl{d\'erivation}
\textsl{int\'erieure} de $\mathfrak{g}_m$. }

\bigskip
	\textbf{2)} \textit{$\hbox{Int}(\mathfrak{g}_m)$ est la composante neutre de
$\hbox{Int}(\mathfrak{g}_m)$}.

\vskip 7mm
\noindent
\textbf{1.3.5.}~\textsc{Remarque} : Le r\'esultat ci-dessus ne vaut pas seulement pour les alg\`ebres de
Takiff $\mathfrak{g}_m$ construites sur une alg\`ebre de Lie semi-simple $\mathfrak{g}_0$ . Un
examen attentif de la d\'emonstration donn\'ee ci-dessus montre qu'il suffit que $\mathfrak{g}_0$ soit
une alg\`ebre de Lie quadratique ayant la propri\'et\'e de Dixmier, \'enonc\'ee dans 1.3.2.

\vskip 10mm
\begin{large}
\noindent
\textbf{1.4. Espaces riemanniens sym\'etriques}
\end{large}

\bigskip 
\noindent
\textbf{1.4.1.}~Soit $G/K$ un espace riemannien sym\'etrique non compact de rang 1. Le groupe $K$
est compact et op\`ere naturellement dans l'espace tangent au point ``neutre'' de $G/K$, not\'e
$\mathfrak{p}$ comme d'habitude. L'alg\`ebre $\Bbb{R}[\mathfrak{p}]^K$ des fonctions polyn\^omes
$K$-invariantes est $\Bbb{R}[Q]$ où $Q$ est une forme quadratique d\'efinie positive sur
$\mathfrak{p}$. Donc $K'$ est le groupe orthogonal de la forme $Q$. Ci-dessous, on utilise les
notations de \cite{He} et la classification des espaces riemanniens sym\'etriques qui s'y trouve.

\medskip
\begin{enumerate}
\item $SO_0(n,1)/SO(n)$ est l'espace hyperbolique r\'eel de dimension $n$. Ici $K' = O(n)$ et
$K'$ et $K$ ont la m\^eme composante connexe $SO(n)$ (c'est l'exemple 1.1.2. 1. ci-dessus).

\medskip
	\item $SU(n,1)/S(U(n) \times U(1))$ est l'espace hyperbolique hermitien. Ici $K' = O(2n)$ (c'est
pour l'essentiel l'exemple 1.1.2. 2. ci-dessus).

\medskip
\item $Sp(n,1)/Sp(n)\times Sp(1)$ est l'espace hyperbolique quaternionnien (de dimension $4n$).
Donc $K' = O(4n)$.

\medskip
\item $F_4/SO(9)$ est l'espace exceptionnel de rang 1, sa dimension est 16. Donc $K' = O(16)$.
\end{enumerate}

\bigskip
\noindent
\textsc{Conclusion} : $\Bbb{R}[\mathfrak{p}]^K$ caract\'erise $K$ (au sens : $K'$ et $K$ ont la m\^eme
dimension) dans le seul cas de l'espace hyperbolique r\'eel.

\vskip 7mm
\noindent
\textbf{1.4.2.}~D'autres exemples d'espaces riemanniens sym\'etriques apparaissent dans \cite{Ra}
(chapitre V), par exemple $SL(n,\Bbb{R})/SO(n)$, où $K'$ est d\'etermin\'e explicitement, et où il
apparaît que $K$ est caract\'eris\'e par ses polyn\^omes invariants. On verra plus loin, comme application
de la notion \textsl{d'orbite caract\'eristique}, que le groupe $K'$ est enti\`erement
d\'etermin\'e dans les cas :
$$
	SO_0(n,n)/SO(n) \times SO(n), \ SU(n,n)/S(U(n) \times U(n))
$$
avec le m\^eme r\'esultat : $K$ est la composante neutre de $K'$.

\vskip 7mm
\noindent
\textbf{1.4.3.}~Comme dans le cas des alg\`ebres de Lie semi-simples, il y a un r\'esultat g\'en\'eral
concernant les alg\`ebres de Lie sym\'etriques complexes, dû \`a Levasseur et Ushirobira (\cite{L-U}) :
lorsque $G/K$ est un espace sym\'etrique de rang $\geq 2$, le groupe $K$ est caract\'eris\'e par ses
polyn\^omes invariants ($K$ est la composante neutre de $K'$).

\vskip 10mm
\begin{large}
\noindent
\textbf{1.5. Alg\`ebres de Jordan}
\end{large}

\bigskip 
	La liste des alg\`ebres de Jordan simples, de dimension finie sur $\Bbb{C}$ est la suivante :

\medskip
\noindent
\textbf{(I)}~$J = M_n(\Bbb{C})$, avec la multiplication de Jordan :
$$
	\{ x,y\} = {1\over 2}\, (xy + yx)
$$
\medskip
\noindent
\textbf{(II)}~$J = MS_n(\Bbb{C})$ est la sous-alg\`ebre de Jordan de $M_n(\Bbb{C})$ constitu\'ee par les
matrices sym\'etri\-ques.

\medskip
\noindent
\textbf{(III)}~$J$ est la sous-alg\`ebre de Jordan de $M_{2n}(\Bbb{C})$ constitu\'ee par les matrices $x$
telles que $x = E\, {}^tx\, E^{-1}$, avec : 
$$
	E = \begin{pmatrix} 
0		&I_n\\
-I_n		&0\\
\end{pmatrix}
$$

\medskip
\noindent
\textbf{(IV)}~$J$ est l'alg\`ebre construite sur l'espace vectoriel $\Bbb{C}^n \otimes \Bbb{C}.1$ où 1
est l'\'el\'ement neutre, et lorsque $x = (x_1,\ldots , x_n)$, $y = (y_1,\ldots ,y_n)$ sont dans
$\Bbb{C}^n$ : 
$$
	xy = \, <x,y>1, \ <x,y> = x_1y_1 +\cdots + x_n y_n
$$

\medskip
\noindent
\textbf{(V)}~$J$ est l'alg\`ebre exceptionnelle constitu\'ee par les matrices hermitiennes \`a 3 lignes et 3
colonnes, \`a coefficients dans l'alg\`ebre de Cayley.

\bigskip
\begin{enumerate}
\item Dans le cas (I), il est connu que le groupe $\hbox{Aut}(J)$ des automorphismes de $J$ est
engendr\'e par le groupe des automorphismes (int\'erieurs) de l'alg\`ebre associative $M_n(\Bbb{C})$
\textsl{et} la transposition des matrices. Dans les cas (II) et (III), il est aussi connu que tout
automorphisme de Jordan de $J$ s'\'etend en un automorphisme de Jordan de $M_p(\Bbb{C})$
($p=n$ dans le cas (II), $p=2n$ dans le cas (III)). Par cons\'equent, l'alg\`ebre
$\Bbb{C}[J]^{\hbox{Aut}(J)}$ des fonctions polyn\^omes $\hbox{Aut}(J)$-invariantes sur $J$, est celle
engendr\'ee par les coefficients du polyn\^ome caract\'eristique. En sens inverse, il est prouv\'e dans
\cite{Ra} (chapitre V, \S 3) que dans tous les cas (I), (II), (III), le groupe $\hbox{Aut}(J)$ est
caract\'eris\'e par ses polyn\^omes invariants. On notera que les 2 polyn\^omes $tr(x^2)$ et $tr(x^3)$
suffisent \`a caract\'eriser $\hbox{Aut}(J)$.

\medskip
\item Dans le cas (IV), le groupe $\hbox{Aut}(J)$ s'identifie au groupe orthogonal $O(n,\Bbb{C})$
agissant sur les $n$ premi\`eres variables $x_1,\ldots , x_n$ de $\Bbb{C}^n$ et trivialement sur
$\Bbb{C}.1$, et l'alg\`ebre des polyn\^omes invariants est engendr\'ee par la derni\`ere coordonn\'ee
$\lambda(x+ \alpha .1)= \alpha$ et par la forme quadratique $Q(x_1,\ldots , x_n, \alpha) = x^2_1
+ \cdots + x^2_n$. On v\'erifie alors que $\hbox{Aut}(J)$ est le groupe des automorphismes de
l'espace vectoriel complexe $J$ qui laisse invariant $\lambda$ et $Q$.

\medskip
\item Dans le cas (V), il est prouv\'e par Chevalley et Schafer (\cite{Ch-Sch}) que l'alg\`ebre de Lie des
d\'erivations de $J$ (de type $F_4$) est caract\'eris\'ee par les 2 polyn\^omes $\hbox{Aut}(J)$-invariants
$tr(xox)$ et $tr(xoxox)$ (où $o$ est la multiplication de Jordan). On a donc la :
\end{enumerate}

\bigskip
\noindent
\textbf{Proposition} : \textit{Soit $J$ une alg\`ebre de Jordan simple, de dimension finie sur
$\Bbb{C}$. Le groupe $\hbox{Aut}(J)$ est caract\'eris\'e par ses polyn\^omes invariants.}

\bigskip
\noindent
\textsc{Remarque} : Cette situation des alg\`ebres de Jordan n'est pas sans relation avec celle des
espaces sym\'etriques : \`a un espace hermitien sym\'etrique correspond un STJ (syst\`eme triple de
Jordan) hermitien positif, et vice-versa, et dans le cas des espaces hermitiens de type tube, ce STJ
``est'' une alg\`ebre de Jordan.

\medskip
	On est donc amen\'e \`a poser la :

\bigskip
\noindent
\textsc{Question} : Le groupe des automorphismes d'un STJ est-il caract\'eris\'e par ses invariants
poly\-n\^omes ?

\vskip 12mm

\section{Sur les orbites caract\'eristiques}

	Soit $G$ un sous-groupe de $GL(V)$. Lorsque $\Omega$ est un \textsl{sous-ensemble} de $V$,
on notera $St(\Omega)$ (le ``stabilisateur'' de $\Omega$) le sous-groupe de $GL(V)$ constitu\'e par
les $f$ dans $GL(V)$ tels que $f(\Omega) = \Omega$. En particulier, lorsque $\Omega$ est une
$G$-orbite, on a \'evidemment : $G \subset St(\Omega)$, et on dira que $\Omega$ est une
\textsl{orbite caract\'eristique} de $G$ si $G = St(\Omega)$ ou au moins si $G$ et $St(\Omega)$
ont m\^eme composante neutre. On notera que dans le cas ``r\'eel'', si $G=K$ est un sous-groupe
\textsl{compact} de $GL(V)$, on a : $K \subset K' \subset St(\Omega)$ pour toute $K$-orbite
$\Omega$. La d\'etermination de $St(\Omega)$ est donc utile pour connaître $K'$, et en particulier :
s'il existe une orbite caract\'eristique, alors $K$ est caract\'eris\'e par ses polyn\^omes invariants.

\vskip 10mm
\begin{large}
\noindent
\textbf{2.1. Orbites caract\'eristiques dans les espaces sym\'etriques}
\end{large}

\bigskip
\noindent
\textbf{2.1.1.}~On consid\`ere l'espace riemannien sym\'etrique $SO_0(n,n)/SO(n) \times SO(n)$.
L'espace $\mathfrak{p}$ s'identifie \`a l'espace $M_n(\Bbb{R})$ des matrices $n \times n$, \`a
coefficients r\'eels, et un couple $(x,y)$ dans $SO(n) \times SO(n)$ op\`ere dans $M_n(\Bbb{R})$ de la
mani\`ere suivante : $(x,y).X = xXy^{-1}$ pour tout $X$ dans $M_n(\Bbb{R})$. Notons $\Omega \ (=
SO(n))$ la $K$-orbite de la matrice unit\'e $I_n$. La d\'etermination de $St(\Omega)$ a \'et\'e faite dans
\cite{Dj}. Il apparaît que si $f(SO(n)) = SO(n)$, alors $f$ est de l'une des deux formes suivantes : 

	ou bien $f(X) = uXv$ avec $u$ et $v$ dans $SO(n)$

	ou bien $f(X) = u\ {}^t\!X v$ avec $u$ et $v$ dans $SO(n)$

\medskip
\noindent
(où ${}^t\!X$ est la transpos\'ee de la matrice $X$).

\medskip
	Dans cet exemple, on voit que $St(\Omega)$, avec $\Omega = SO(n)$, admet $K$ comme
sous-groupe normal connexe d'indice 2.

\vskip 7mm
\noindent
\textbf{2.1.2.}~Consid\'erons l'espace sym\'etrique $SU(n,n)/S(U(n)\times U(n))$ ; l'espace
$\mathfrak{p}$ \'etant identifi\'e \`a $M_n(\Bbb{C})$, notons $\Omega = U(n)$ l'orbite sous $K$ de la
matrice unit\'e $I_n$. Dans \cite{M}, il est prouv\'e que les $f : \mathfrak{p} \longrightarrow
\mathfrak{p}$, \textsl{$\Bbb{C}$-lin\'eaires} et bijectives telles que $f(U(n)) = U(n)$ sont celles
ayant l'une des deux formes suivantes :

\medskip
	\quad \ \  $f(x) = uxv^*$ avec $u$ et $v$ unitaires

\medskip
	ou \ $f(x) = u{}^t\!x v^*$ avec $u$ et $v$ unitaires.

\bigskip
	Rappelons que si $u$ et $v$ sont 2 matrices unitaires, il existe 2 autres matrices unitaires $u_0$ et
$v_0$ telles que $det(u_0 v_0) = 1$ et $uxv^* = u_0xv^*_0$ pour tout $x$ ; par suite, les bijections
lin\'eaires complexes $f$ ayant la premi\`ere forme ci-dessus, sont exactement celles provenant de
l'action naturelle de $K$ dans $\mathfrak{p}$. Si on note $St_{\Bbb{C}}(\Omega)$ le sous-groupe de
$St(\Omega)$ constitu\'e par les bijections $\Bbb{C}$-lin\'eaires, ce qui pr\'ec\`ede montre que $K$ (en
fait l'image de $K$ dans $GL(\mathfrak{p})$) est la composante neutre de $St_{\Bbb{C}}(\Omega)$.
Mais qu'en est-il du groupe $St(\Omega)$ ? Un compl\'ement au th\'eor\`eme de Marcus cit\'e ci-dessus a
\'et\'e apport\'e par Djokovic dans \cite{Dj} : le groupe $St(\Omega)$ est engendr\'e par
$St_{\Bbb{C}}(\Omega)$ et par la transformation antilin\'eaire $f(x) = \overline{x}$. Donc $K$ est
aussi la composante neutre de $St(\Omega)$ et $\Omega = U(n)$ est une orbite caract\'eristique de $K$.

\vskip 7mm
\noindent
\textbf{2.1.3.}~Soit $G/K$ un espace sym\'etrique \textsl{hermitien} (irr\'eductible). L'espace
vectoriel r\'eel $\mathfrak{p}$ admet une structure complexe $K$-invariante, de sorte que l'action de
$K$ dans $\mathfrak{p}$ se fait par des applications lin\'eaires complexes. L'espace homog\`ene $G/K$
se r\'ealise (par le plongement de Harish-Chandra) comme la boule unit\'e $\mathcal{D}$ dans
$\mathfrak{p}_+$, relativement \`a une norme bien pr\'ecise (voir par exemple \cite{Sa}, chapitre II, \S
4). Dans la fronti\`ere de $\mathcal{D}$ se trouve une $K$-orbite particuli\`ere qui est la
\textsl{fronti\`ere de Shilov} du domaine born\'e sym\'etrique $\mathcal{D}$. Notons $\Omega$
cette orbite particuli\`ere. Dans l'exemple pr\'ec\'edent, $\mathcal{D}$ est l'ensemble des matrices $x$
dans $M_n(\Bbb{C})$ telles que $(I_n - xx^*)$ soit \textsl{d\'efinie positive}, la norme est la
norme spectrale, et $\Omega$ est l'ensemble $U(n)$ des matrices unitaires. Dans le cas g\'en\'eral, il se
trouve que la fronti\`ere de Shilov $\Omega$ de $\mathcal{D}$ est l'ensemble des points extr\'emaux de
$\overline{\mathcal{D}}$. D\`es lors, soit $f$ un automorphisme $\Bbb{R}$-lin\'eaire de $\mathfrak{p}_+$
telle que : $f(\Omega) = \Omega$ ; on a alors $f(\overline{\mathcal{D}}) \subset
\overline{\mathcal{D}}$, d'où : $\| f(x)\| \leq \|x\|$ pour tout $x$ dans $\mathfrak{p}_+$ (où $\| \ 
\|$ est la norme d\'efinissant la boule unit\'e $\mathcal{D}$). Comme la bijection r\'eciproque $g =
f^{-1}$ de $f$ a la m\^eme propri\'et\'e : $g(\Omega) = \Omega$, on a : $\| g(x) \| \leq \| x \|$, d'où
$\|x\| \leq \|f(x)\|$ et ainsi $\|f(x)\| = \|x\|$ pour tout $x$ dans $\mathfrak{p}_+$. Donc
$f(\mathcal{D}) = \mathcal{D} = g(\mathcal{D})$. \textsl{Supposons dor\'enavant que $f$ soit
$\Bbb{C}$-lin\'eaire}, de sorte que $f$ est un automorphisme biholomorphe de $\mathcal{D}$. La
comparaison du groupe $\hbox{Hol}(\mathcal{D})$ des automorphismes biholomorphes de
$\mathcal{D}$ avec le groupe $I(\mathcal{D})$ des isom\'etries (relativement \`a la m\'etrique de
Bergman) de
$\mathcal{D}$ est \'etablie dans \cite{Sa}, page 88 : $\hbox{Hol}(\mathcal{D})$ est un sous-groupe
d'indice 2 de $I(\mathcal{D})$, et si on note $I(\mathcal{D})^0$ la composante neutre de
$I(\mathcal{D})$, on a : 
$$
	I(\mathcal{D})^0 \subset \hbox{Hol}(\mathcal{D}) \subset I(\mathcal{D})
$$
et $I(\mathcal{D})^0$ est un sous-groupe d'indice au plus 2 de $\hbox{Hol}(\mathcal{D})$. Lorsque
$f \in I(\mathcal{D})^0$, ce qui est le cas lorsque $\hbox{Hol}(\mathcal{D}) = I(\mathcal{D})^0$, on
a imm\'ediatement : $f \in G$ et $f(0) = 0$, d'où : $f \in K$. Dans tous les cas,
$St_{\Bbb{C}}(\Omega)^0 \subset K$ et comme $K \subset ST_{\Bbb{C}}(\Omega)^0$, il apparaît
que $K$ est la composante neutre de $St_{\Bbb{C}}(\Omega)$. 

	Ceci ne r\'epond pas \`a la question de savoir si $\Omega$ est une orbite caract\'eristique de $K$,
comme c'est le cas dans 2.1.2. On est donc amen\'e \`a faire lorsque le rang est $\geq 2$, la 

\bigskip
\noindent
\textsc{Conjecture} : $K$ est la composante neutre de $St(\Omega)$.

\vskip 7mm
\noindent
\textbf{2.1.4.}~Compte-tenu de l'exemple trait\'e dans 2.1.1., on peut poser la :

\bigskip
\noindent
\textsc{Question} : Lorsque $G/K$ est un espace riemannien sym\'etrique de rang $\geq 2$, le groupe
$K$ admet-il une orbite caract\'eristique dans $\mathfrak{p}$ ? Si oui, quelles sont les orbites
caract\'eristiques ?

\vskip 10mm
\begin{large}
\noindent
\textbf{2.2. Orbites adjointes caract\'eristiques dans les alg\`ebres de Lie semi-simples}
\end{large}

\bigskip
\noindent
\textbf{2.2.1.}~Soient $\mathfrak{g} = \mathfrak{g}l(n,\Bbb{C})$, $\Omega$ l'orbite adjointe d'un
\'el\'ement g\'en\'erique fix\'e de $\mathfrak{g}$ et $\mathcal{N}$ le c\^one des \'el\'ements nilpotents de
$\mathfrak{g}$. Dans \cite{Wa}, il est prouv\'e que $St(\Omega)$ est engendr\'e par
$\hbox{Int}(\mathfrak{g})$ et par la transposition des matrices. Dans \cite{B-P-Wa}, il est prouv\'e que
$St(\mathcal{N})$ est engendr\'e par $\hbox{Int}(\mathfrak{g})$, la transposition des matrices
\textsl{et} l'ensemble des homoth\'eties $x \mapsto \lambda x$ avec $\lambda$ dans
$\Bbb{C}^*$.

	Une sorte de r\'ecapitulatif r\'ecent sur ces r\'esultats se trouve dans \cite{Pl-Dj}.

\vskip 7mm
\noindent
\textbf{2.2.2.}~Soit $\mathfrak{g}$ une alg\`ebre de Lie semi-simple complexe, de dimension finie, de
rang $r$. On note $(p_1,\ldots ,p_r)$ un syst\`eme de g\'en\'erateurs homog\`enes, alg\'ebriquement
ind\'ependants de l'alg\`ebre $\Bbb{C}[\mathfrak{g}]^{\hbox{Int}(\mathfrak{g})}$, et $\pi :
\mathfrak{g} \longrightarrow \Bbb{C}^r$ l'application d\'efinie par : $\pi(x) = (p_1(x),\ldots 
p_r(x))$ pour tout $x$.

\bigskip
	Soit maintenant $f$ dans $GL(\mathfrak{g})$ tel que $pof = p$ pour tout polyn\^ome
$\hbox{Int}(\mathfrak{g})$-invariant sur $\mathfrak{g}$ (autrement  dit : $f \in
\hbox{Int}(\mathfrak{g})'$). On a imm\'ediatement : 
$$
	{}^t\!fodpof = dp\quad \hbox{pour tout} \ p \ \hbox{dans}\ \Bbb{C}[p_1,\ldots ,p_r]
$$
Notons $d(x)$ le rang de $\pi$ au point $x$ (pour $x$ dans $\mathfrak{g}$). Il vient alors :
$d(f(x)) = d(x)$ pour tout $x$, et $f(\mathfrak{g}_r) = \mathfrak{g}_r$, où $\mathfrak{g}_r$ est
l'ensemble des \'el\'ements r\'eguliers de $\mathfrak{g}$. De m\^eme un tel $f$ transforme un \'el\'ement
\textsl{g\'en\'erique} en un \'el\'ement g\'en\'erique, car $f$ laisse invariant le polyn\^ome discriminant de
$\mathfrak{g}$. D'apr\`es les th\'eor\`emes de Kostant, il vient : $\hbox{Int}(\mathfrak{g})' \subset
St(\Omega)$ pour toute orbite r\'eguli\`ere $\Omega$ et de m\^eme pour toute orbite g\'en\'erique.

\vskip 7mm
\noindent
\textbf{2.2.3.}~Soit $\mathcal{N}_0$ l'orbite r\'eguli\`ere nilpotente. Comme $\mathcal{N}_0$ est
dense dans $\mathcal{N}$, on a : $St(\mathcal{N}_0) \subset St(\mathcal{N})$.

\bigskip
\noindent
\textbf{Lemme} : \textit{Soit $f$ dans $GL(\mathfrak{g})$. Sont \'equivalentes :}

	(1) \quad $f(\mathcal{N}_0)= \mathcal{N}_0$

\medskip
	(2) \quad \textit{$f(\mathcal{N})= \mathcal{N}$ et $d(f(x)) = d(x)$ pour tout $x$ dans
$\mathcal{N}$.}

\bigskip
\noindent
\textsc{D\'emonstration} : \textbf{(i)} Soit $f$ tel que $f(\mathcal{N}_0) = \mathcal{N}_0$ ; on a donc
: $f(\mathcal{N}) = \mathcal{N}$. Soit $p$ un polyn\^ome $\hbox{Int}(\mathfrak{g})$-invariant, sans
terme constant ; alors $pof$ est nul sur $\mathcal{N}$ et il existe des fonctions polyn\^omes
$\varphi_1,\ldots , \varphi_r$ telles que $pof = \sum\, \varphi_j\, p_j$, d'où : 
$$
	{}^t\!f o dp o f = \sum\, \varphi_j\, dp_j + \sum\, p_j\, d\varphi_j
$$
et, pour tout $x$ dans $\mathcal{N}$ : 
$$
	{}^t\!f(dp(f(x)) = \sum\, \varphi_j(x)\, dp_j(x)
$$
soit :
$$
	dp(f(x)) = \sum\, \varphi_j(x) \check{f}(dp_j(x))\quad \hbox{avec}\quad \check{f} =
({}^t\!f)^{-1}.
$$
Notons $E(x)$ le sous-espace vectoriel de $\mathfrak{g}^*$ engendr\'e par les formes lin\'eaires
$dp_1(x),\ldots , dp_r(x)$. Ce qui pr\'ec\`ede montre que : $E(f(x)) \subset \check{f}(E(x))$ et $d(f(x))
\leq d(x)$. Mais on a aussi, avec $g = f^{-1} : E(g(y)) \subset \check{g}(E(y))$, ce qui s'\'ecrit, avec $y
= f(x) : E(x) \subset \check{g} (E(f(x)))$, d'où $d(x) \leq d(f(x))$ et ainsi :
$$
	\check{f}\ \hbox{induit un isomorphisme d'espaces vectoriels}
$$
de $E(x)$ sur $E(f(x)),\ d(f(x)) = d(x)$, pour tout $x$ dans $\mathcal{N}$. 

\vskip 4mm
\textbf{(ii)} Sous l'hypoth\`ese (2), on a : $f(\mathcal{N}_0) \subset \mathcal{N}_0$ (d'où
$f(\mathcal{N}_0) = \mathcal{N}_0$) car $\mathcal{N}_0$ est l'ensemble des $x$ dans
$\mathcal{N}$ tels que : $d(x) = r$.

\vskip 5mm
\noindent
\textsc{Conjectures} : \textbf{(1)} Soit $\Omega$ une orbite g\'en\'erique. Alors : 
$$
	St(\Omega)^0 = \hbox{Int}(\mathfrak{g})
$$
\textbf{(2)} Soient $\mathcal{N}$ le c\^one des \'el\'ements nilpotents de $\mathfrak{g}$ et
$\mathcal{N}_0$ l'orbite nilpotente r\'eguli\`ere. Alors :
$$
	St(\mathcal{N})^0 = St(\mathcal{N}_0)^0 = \hbox{Int}(\mathfrak{g}) \times \Bbb{C}^*.
$$

\vfill\eject

\section{Sur les distributions invariantes caract\'eristiques}

\begin{large}
\noindent
\textbf{3.1. Alg\`ebres de Lie nilpotentes}
\end{large}

	Soit $\mathfrak{g}$ une alg\`ebre de Lie nilpotente, de dimension finie sur $\Bbb{R}$, et soit $G$ un
groupe de Lie connexe d'alg\`ebre de Lie $\mathfrak{g}$. Le groupe $G$ op\`ere dans $\mathfrak{g}$
(resp.  $\mathfrak{g}^*$) au moyen de la repr\'esentation adjointe (resp. coadjointe), et on notera
$Ad(G)$, (resp. $Ad^*(G)$) l'image de $G$ dans $GL(\mathfrak{g})$ (resp. $GL(\mathfrak{g}^*)$)
r\'esultant de la repr\'esentation adjointe (resp. coadjointe). On examine ici les questions pos\'ees
ci-dessus pour la repr\'esentation coadjointe de $G$ dans $\mathfrak{g}^*$.

\vskip 7mm
\noindent
\textbf{3.1.1.}~On a vu dans l'exemple 1.1.3. que le groupe de Heisenberg n'est pas caract\'eris\'e par
ses polyn\^omes invariants. On rappelle qu'ici $\mathfrak{g}^*$ s'identifie \`a $\Bbb{R}^3$, avec les
coordonn\'ees $(x,y,z)$, et que les champs de vecteurs coadjoints sont les combinaisons lin\'eaires de
$$
	A = z\, {\partial\over \partial x}, \quad B = z\, {\partial\over \partial y}
$$
Les \textsl{distributions} $T$, invariantes sous l'action coadjointe de $G$ sur $\mathfrak{g}^*$,
s'\'ecrivent : 
$$
	T = 1_{x,y} \otimes U(z) + V(x,y) \otimes \delta(z)
$$
où $U$ et $V$ sont des distributions arbitraires. Parmi ces distributions invariantes, on privil\'egie les
distributions de type :
$$
	T_1 = 1_{x,y}\ \otimes \delta(z-z_0)\quad (z_0 \not= 0)
$$
$$
	T_2 = \delta(x-x_0) \otimes \delta(y-y_0) \otimes \delta(z) \quad (x_0, y_0) \in \Bbb{R}^2
$$
qui sont des \textsl{mesures} invariantes port\'ees par les orbites coadjointes de $G$ dans
$\mathfrak{g}^*$. 

	Soit $ L = a\, {\partial\over \partial x} + b\, {\partial\over \partial y} + c\, {\partial\over
\partial z}$ un champ de vecteurs sur $\mathfrak{g}^*$, \`a coefficients dans $C^\infty(\Bbb{R}^3)$, qui
annule les mesures invariantes port\'ees par les $G$-orbites.

\vskip 4mm
\noindent
- \begin{eqnarray*}
0 &=	&L(1_{x,y} \otimes \delta(z-z_0)) = 1_{x,y} \otimes c(x,y,z) \delta '(z_0)\\
		&=	&1_{x,y} \otimes [c(x,y,z_0) \delta '(z_0) - {\partial\over \partial z} c(x,y,z_0) \delta(z_0)]
\end{eqnarray*}
On a donc, pour tout couple $(\psi, \theta)$ de fonctions $C^\infty$ (\`a support compact) $\psi =
\psi(x,y), \theta = \theta(z)$ :
$$
	0 = \int_{\Bbb{R}^2}\, \psi (x,y) [- c(x,y,z_0) \theta '(z_0) - {\partial c\over \partial z} (x,y,z_0)
\theta(z_0)]dx dy
$$
En choisissant $\theta$ telle que : $\theta(z_0) = 0$ et $\theta '(z_0) = 1$, il vient :
$$
	\int \psi(x,y) c(x,y,z_0)dx dy = 0\ \hbox{pour toute}\ \psi \ \hbox{dans}\ \mathcal{D}(\Bbb{R}^2).
$$
Donc : $c(x,y,z_0) = 0$ pour tout $(x,y,z_0)$ dans $\Bbb{R}^2 \times \Bbb{R}^*$, et ainsi $c = 0$.

\vskip 7mm
\noindent
- \quad $0 = L(\delta(x-x_0) \otimes \delta(y-y_0) \otimes \delta(z)) = (L_0\, \delta(x-x_0) \otimes
\delta(y-y_0)) \otimes \delta(z)$

\medskip
\noindent
avec $L_0 = a(x,y,0) {\partial\over \partial x} + b(x,y,0) {\partial\over \partial y}$.

\medskip
	Donc $L_0(\delta(x-x_0)\otimes \delta(y-y_0)) = 0$ pour tout $(x_0, y_0)$ dans
$\Bbb{R}^2$, et l'op\'erateur diff\'erentiel ${}^t\! L_0$ annule toutes les fonctions $\psi $ $C^\infty$ sur
$\Bbb{R}^2$ :
$$
	<{}^t\! L_0(\psi),\ \delta(x-x_0) \otimes \delta(y-y_0)>\,  = 0.
$$
Donc ${}^t\! L_0 = - \big ( {\partial a\over \partial x} (x,y,0) + {\partial b\over \partial y}
(x,y,0)\big ) - a {\partial\over \partial x} - b {\partial \over \partial y} = 0$ et par suite : $a(x,y,0) =
b(x,y,0) = 0$ pour tous $x$ et $y$, et il existe 2 fonctions $\alpha$ et $\beta$, de classe $C^\infty$,
telles que :
$$
	a(x,y,z) = z \alpha(x,y,z),\ b(x,y,z) = z \beta(x,y,z).
$$
Donc $L = \alpha A + \beta B$.

\vskip 7mm
\noindent
\textbf{3.1.2.}~La question de savoir si on a le m\^eme r\'esultat pour une alg\`ebre de Lie nilpotente
quelconque a \'et\'e trait\'ee par J. Dixmier (\cite{Di-2}, \S 4), qui montre que la r\'eponse est n\'egative, en
exhibant le contre-exemple suivant : 

\bigskip
	- Soit $\mathfrak{g}$ l'alg\`ebre de Lie (nilpotente), de dimension 6 sur $\Bbb{R}$, admettant la base
$(e_1,\ldots ,e_6)$ avec :
$$
	[e_1,e_2] = e_5, \quad [e_1,e_3] = e_6, \quad [e_2, e_4] = e_6.
$$
\indent
On note $(y_1,\ldots , y_6)$ les fonctions coordonn\'ees sur $\mathfrak{g}^*$, correspondant \`a la
base duale de $(e_1,\ldots , e_6)$. Les champs coadjoints d\'efinis respectivement par $e_1,\ldots ,
e_6$ sont :
$$
	A_1 = y_5 {\partial\over \partial y_2} + y_6 {\partial\over \partial y_3}, \quad A_2 = -y_5
{\partial\over \partial y_1} + y_6 {\partial\over \partial y_4}, \quad A_3 = -y_6 {\partial\over
\partial y_1}, \quad
$$

$$
	A_4 = -y_6 {\partial\over \partial y_2}, \quad A_5 = A_6 = 0
$$
Dixmier montre que le champ $y_5 {\partial\over \partial y_1}$ annule toutes les mesures
invariantes port\'ees par les orbites coadjointes, mais n'est pas combinaison lin\'eaire, \`a coefficients
fonctions $C^\infty$ des champs $A_1,\ldots ,A_6$.

\vskip 4mm
- Parmi les distributions invariantes $T$ sur $\mathfrak{g}^*$, il y a les distributions de l'un des
trois types suivants : 
\begin{eqnarray*}
	T_1	&=		&1(y_1,y_2,y_3,y_4) \otimes U(y_5, y_6)\\
\smallskip \\
	T_2	&=		&y_1 \otimes 1_{y_2} \otimes {\partial V\over \partial y_4} (y_3, y_4, y_5) \otimes
\delta(y_6)\\
\smallskip \\
		&	&-1_{y_1} \otimes y_2 \otimes {\partial V\over \partial y_3} (y_3, y_4, y_5) \otimes \delta(y_6)\\
\smallskip \\
		&	&-1_{y_1,y_2} \otimes y_5  V(y_3, y_4, y_5) \otimes \delta '(y_6)\\
\smallskip \\
	T_3	&=		&1_{y_1,y_2} \otimes W(y_3, y_4,y_5) \otimes \delta(y_6)
\end{eqnarray*}
où $U, V$ et $W$ sont des distributions arbitraires.

\smallskip
	(En fait, toute distribution invariante est une somme $T_1+T_2+T_3$.)

\vskip 4mm
- Soit $L = \sum \varphi_j\, {\partial\over \partial y_j}$ un champ de vecteurs qui annule
toutes les distributions invariantes, les $\varphi_j$ \'etant des fonctions $C^\infty$ (ou des fonctions
polyn\^omes).

\medskip
	\quad . Il annule toutes les distributions de type $T_1$, ce qui s'\'ecrit : $(\varphi_5 \partial_5 +
\varphi_6 \partial_6)U = 0$ pour toute distribution $U$ (on \'ecrira $\partial_j$ \`a la place de
${\partial\over \partial y_j}$). D'où $\varphi_5 = \varphi_6 = 0$.

\medskip
	\quad . Il annule toutes les distributions de type $T_3$, ce qui s'\'ecrit :
$$
	\varphi_3\ 1_{y_1 y_2} \otimes \partial_3 W \otimes \delta(y_6) + \varphi_4\ 1_{y_1 y_2} \otimes
\partial_4 W \otimes \delta(y_6) = 0
$$
d'où, pour toute distribution $W$ : 
$$
	\varphi_3 (y_1,\ldots , y_5,0)1_{y_1y_2} \otimes \partial_3 W + \varphi_4(y_1,\ldots ,
y_5,0)1_{y_1y_2} \otimes \partial_4W = 0
$$
$$
\hbox{et} \quad \varphi_3(y_1,\ldots , y_5,0) = \varphi_4(y_1,\ldots , y_5,0) = 0.
$$
Il existe donc des fonctions ($C^\infty$ ou polyn\^omes suivant le cas) $\psi_3,\ \psi_4$ telles que :
$\varphi_3 = y_6 \psi_3,\ \varphi_4 = y_6 \ \psi_4$, et :
\begin{eqnarray*}
	L &= 	&\varphi_1 \partial_1 + \varphi_2 \partial_2 + \psi_3 y_6 \partial_3 + \psi_4 y_6
\partial_4\\
	&=	&\varphi_1 \partial_1 + \varphi_2 \partial_2 + \psi_3 A_1 + \psi_4 A_2 - \psi_3 y_5 \partial_2
+ \psi_4 y_5 \partial_1\\
	&=	&(\varphi_1 + \psi_4 y_5)\partial_1 + (\varphi_2 - \psi_3 y_5)\partial_2 + \psi_3 A_1 +
\psi_4 A_2
\end{eqnarray*}

\vskip 4mm
	. $L$ annule toutes les distributions de type $T_2$, ce qui revient \`a : $(\theta_1 \partial_1 +
\theta_2 \partial_2)T_2 = 0$ avec $\theta_1 = \varphi_1 + \psi_4 y_5$ et $\theta_2 = \varphi_2 -
\psi_3 y_5$. Donc : 
$$
	\theta_1\, 1_{y_1y_2} \otimes \partial_4 \, V \otimes \delta(y_6) - \theta_2 \, 1_{y_1y_2} \otimes
\partial_3 \ V \otimes \delta(y_6) = 0
$$
pour toute distribution $V$, soit encore :
$$
	\theta_1(y_1,\ldots , y_5, 0)\partial_4 \, V - \theta_2(y_1,\ldots , y_5,0)\partial_3 \, V = 0.
$$
Il existe donc $\psi_1,\, \psi_2$ telles que : $\theta_1 = y_6 \psi_1,\ \theta_2 = y_6 \psi_2$ et :
\begin{eqnarray*}
	L	&=		& \psi_1 y_6 \partial_1 + \psi_2 y_6 \partial_2 + \psi_3 A_1 + \psi_4 A_2\\
	&=		& \psi_3 A_1 + \psi_4 A_2 - \psi_1 A_3 - \psi_2 A_4
\end{eqnarray*}
L'examen d'autres alg\`ebres de Lie nilpotentes am\`ene \`a faire la :

\vskip 4mm
\noindent
\textsc{Conjecture} : Soit $\mathfrak{g}$ une alg\`ebre de Lie nilpotente et soit $L$ un champ de
vecteurs (\`a coefficients $C^\infty$) sur $\mathfrak{g}^*$, qui annule toutes les distributions
invariantes sur $\mathfrak{g}^*$. Alors il existe une fonction $\Phi : \mathfrak{g}^*
\longrightarrow \mathfrak{g}$, de classe $C^\infty$ telle que :
$$
	L(x) = \hbox{ad}^*(\Phi(x)).x \quad \hbox{pour tout} \ x \  \hbox{dans}\ \mathfrak{g}^*.
$$

\vskip 10mm
\begin{large}
\noindent
\textbf{3.2. Groupes unipotents}
\end{large}

\vskip 4mm
\noindent
\textbf{3.2.1.}~Soit $G$ le groupe ``unipotent standard'', c'est-\`a-dire le sous-groupe de
$GL(n,\Bbb{R})$ constitu\'e par les matrices unipotentes triangulaires sup\'erieures $g = (g_{ij})_{i,j}$
avec $g_{ii} = 1 \ (1\leq i \leq n), \ g_{ij} = 0$ lorsque $i>j$.

\smallskip
	On consid\`ere l'op\'eration naturelle de $G$ dans $\Bbb{R}^n$.

	Les champs de vecteurs associ\'es sont les :
$$
	L_{ij} = x_j {\partial\over \partial x_i}\quad 1 \leq i \leq n-1, \quad i<j\leq n
$$
c'est-\`a-dire : 
\begin{eqnarray*}
	x_2\, {\partial\over \partial x_1},\ldots , x_n\, {\partial\over \partial x_1}\\
	x_3\, {\partial\over \partial x_2},\ldots , x_n\, {\partial\over \partial x_2}\\
	\ldots \ldots \ldots \ldots \ldots \ldots\\
		x_n \, {\partial\over \partial x_{n-1}}
\end{eqnarray*}
	Les fonctions polyn\^omes en $x_1,\ldots ,x_n$, qui sont $G$-invariantes, sont celles qui ne
d\'ependent que de $x_n$ ; c'est l'alg\`ebre $\Bbb{R}[x_n]$, et le groupe $G'$ est celui des matrices de
la forme : 
$$
	\begin{pmatrix}
	*	&	&\cdots	&*\\
		&	&	&\\
	*	&*	&\cdots	&*\\
	0	&0	&\cdots	&1\\
\end{pmatrix}
$$

\bigskip
\noindent
\textsl{D\'etermination des distributions invariantes $T$}

\medskip
\noindent
$\bullet$ \begin{eqnarray*}
	x_2		& \displaystyle {\partial T\over \partial x_1} &=	\ldots = x_n {\partial T\over \partial x_1} = 0\\
\smallskip \nonumber\\
								&\displaystyle {\partial T\over \partial x_1} &= U_0(x_1) \otimes \delta(x_2,\ldots ,x_n)\\
\smallskip \nonumber\\
							&T	&= U(x_1) \otimes \delta(x_2,\ldots , x_n) + 1_{x_1} \otimes V(x_2, \ldots ,x_n)
\end{eqnarray*}

\medskip
\noindent
$\bullet$ $T_1 = U(x_1) \otimes \delta (x_2,\ldots ,x_n)$ est invariante car
$$
	x_j {\partial\over \partial x_i} \delta(x_2,\ldots ,x_n) = 0
$$
lorsque $2 \leq i \leq n-1$ et $i < j \leq n$

\bigskip
\noindent
$\bullet$ $T_2 = 1_{x_1} \otimes V(x_2,\ldots ,x_n)$ est annul\'ee par les $L_{ij}$ et en particulier :
$$
	0 = 1_{x_1} \otimes L_{ij} V \quad \hbox{avec :} \quad 2 \leq i \leq n-1 \quad \hbox{et} \quad
i<j\leq n.
$$
Donc $V$ est invariante par le sous-groupe unipotent standard de $GL(n-1, \Bbb{R})$ op\'erant dans
$\Bbb{R}^{n-1} = \{(0, x_2,\ldots ,x_n)|(x_2,\ldots ,x_n) \in \Bbb{R}^{n-1}\}$.

\vskip 7mm
\noindent
\textbf{Proposition} : \textit{Soit $n\geq 2$. Soit $L$ un champ de vecteurs \`a param\`etres sur
$\Bbb{R}^n$ qui annule toutes les distributions invariantes. Alors :
$$
	L = \sum \varphi_{ij}\, L_{ij}
$$}

\vskip 4mm
\noindent
\textsc{D\'emonstration} : On traite le cas $n=2$. Les distributions invariantes s'\'ecrivent : 
$$
	T = U(x_1) \otimes \delta(x_2) + 1_{x_1} \otimes V(x_2)
$$
où $U$ et $V$ sont des distributions arbitraires.

\vskip 5mm
	Soit $L = a(w,x_1,x_2){\partial\over \partial x_1} + b(w,x_1,x_2){\partial\over \partial
x_2}$ ($w$ est le param\`etre).
 
\bigskip
$\bullet$ $L(1_{x_1} \otimes V(x_2)) = 0$, avec $V$ arbitraire (par exemple $V = x_2)$, donne
$b=0$.

\bigskip
$\bullet$ $a(w,x_1,x_2) {dU\over dx_1} \otimes \delta(x_2) = 0$ avec $U$ arbitraire
\'equivaut \`a :
$$
	a(w,x_1,0) {dU\over dx_1} = 0, \quad \hbox{avec} \quad U \quad \hbox{arbitraire,}
$$
d'où $a(w,x_1,0)=0$, et il existe une fonction $\varphi(w,x_1,x_2)$ telle que : $a(w,x_1,x_2) = x_2
\varphi(w,x_1,x_2)$.

\smallskip
	Donc $L = \varphi\, x_2\, {\partial\over \partial x_1}$.

\vskip 5mm
\noindent
$\bullet$ On traite le cas g\'en\'eral, par r\'ecurrence sur $n$.
$$
	L = a_1(w,x){\partial\over \partial x_1} + a_2(w,x) {\partial\over \partial x_2} +\cdots + a_n(w,x)
{\partial\over \partial x_n}
$$
$$
	T = U(x_1) \otimes \delta(x_2,\ldots ,x_n) + 1_{x_1} \otimes V(x_2,\ldots ,x_n)
$$
Par l'hypoth\`ese de r\'ecurrence :
$$
	L_0 = \sum^n_{j=2} \, a_j(w,x){\partial\over \partial x_j}
$$
s'\'ecrit : 
$$
	 L_0 = \sum\, \varphi_{ij}(w,x)L_{ij}(x) \quad (2\leq i).
$$

\smallskip
D'où :
$$
	a_1 {\partial\over \partial x_1}. U(x_1) \otimes \delta(x_2,\ldots ,x_n) = 0
$$
et :
$$
	a_1(w,x_1,0,\ldots ,0){dU\over dx_1} = 0
$$
pour toute distribution $U$. Il existe donc des fonctions $\varphi_2(w,x),\ldots , \varphi_n(w,x)$
telles que :
$$
	a_1(w,x) = x_2\, \varphi_2(w,x) + x_3\, \varphi_3(w,x)+\cdots + x_n \varphi_n(w,x).
$$
D'où : $L = \varphi_2 x_2 {\partial\over \partial x_1} +\cdots + \varphi_n x_n
{\partial\over
\partial x_1} + \sum_{i\geq 2} \varphi_{ij} (w,x)L_{ij}(x)$.

\vfill\eject
\noindent
\textbf{3.2.2.}~Ici $G$ est le sous-groupe unipotent de $GL(n,\Bbb{R})$ :
$$
	G = \{exp(tX)\, | \, t \in \Bbb{R}\}
$$
$$
	X = \begin{pmatrix}
	0	&1	&0	&\ldots 	&0\\
	0	&0	&1	& 			&\vdots\\
	\vdots	&\vdots	&\ddots	&\ddots	&0\\
	\vdots	&\vdots	&\ddots	&\ddots	&1\\
	0	&0	&0	&\ldots	&0\\
\end{pmatrix}
$$
(Si $(e_1,\ldots ,e_n)$ est la base canonique de $\Bbb{R}^n$, l'endomorphisme $f$ de $\Bbb{R}^n$
dont la matrice est $X$, v\'erifie : 
$$
	f(e_1) = 0\quad f(e_j) = e_{j-1} \quad (2\leq j \leq n)), \quad (n\geq 3).
$$

\medskip
	La repr\'esentation naturelle de $G$ dans $\Bbb{R}^n$ est la restriction \`a 
$N = \Big\{ \begin{pmatrix} 1	&t\\ 0	&1\\ \end{pmatrix}\ | \ t\in \Bbb{R}\Big \}$ de la
repr\'esentation irr\'eductible de dimension $n$ de $SL(2,\Bbb{R})$.

\bigskip
\begin{enumerate}
\item Soit $\Omega$ l'ensemble des $x = (x_1,\ldots ,x_n)$ tels que $x_n \not = 0$ et soit
$\mathcal{H}$ l'hyperplan d'\'equation $x_{n-1} = 0$. Pour chaque $x$ dans $\Omega$, la $G$-orbite
$G.x$ rencontre $\mathcal{H}$ en un point et un seul $R(x)$, dont les coordonn\'ees dans la base
canonique : $Q_1,\ldots , Q_{n-2}, Q_{n-1}$ (on a donc :
$$
	G.x \cap \mathcal{H} = R(x) = (Q_1(x), Q_2(x),\ldots ,Q_{n-2}(x), 0, Q_{n-1}(x))
$$
sont des fonctions rationnelles (de $x$ dans $\Omega$), qui s\'eparent les orbites de $G$ dans
$\Omega$, et en fait engendrent le corps des fonctions rationnelles $G$-invariantes sur $\Bbb{R}^n.$

\bigskip
\item Pour chaque $x$ dans $\Bbb{R}^n$, on pose :
\begin{eqnarray*}
	P_r(x) 	&=		&x^{n-r-1}_n\, Q_r(x) \quad 1\leq r \leq n-2\\
P_{n-1}(x)	&=	&x_n\\
\end{eqnarray*}

Les fonctions $P_1,\ldots ,P_{n-1}$ sont en fait des fonctions polyn\^omes $G$-invariantes,
alg\'ebrique\-ment ind\'ependantes qui engendrent le corps des fonctions rationnelles invariantes, mais
ne constituent pas un syst\`eme de g\'en\'erateurs de l'alg\`ebre des polyn\^omes invariants. Pr\'ecis\'e\-ment, si
$P$ est une fonction polyn\^ome invariante, il existe un entier naturel $k$ tel que :
$$
	x^k_n P \in \Bbb{R}[P_1,\ldots ,P_{n-1}].
$$
Les polyn\^omes $P_1,\ldots ,P_{n-1}$ peuvent \^etre explicitement \'ecrits dans les variables $x_1,\ldots
,x_n$, par exemple, lorsque $n=3$ :
$$
	P_1(x_1,x_2,x_3) = x_1x_3 - {1\over 2} x^2_2
$$
$$
	P_2(x_1,x_2,x_3) = x_3
$$

\bigskip
\item On peut alors d\'emontrer, par r\'ecurrence sur $n\geq 3$ le r\'esultat suivant : Soit $L =
a_1 {\partial\over \partial x_1} +\cdots + a_n {\partial\over \partial x_n}$ un champ de vecteurs sur
$\Bbb{R}^n$ qui annule les fonctions $P_1,\ldots ,P_n$. Il existe alors une fonction $\varphi$ telle
que : 
$$
	L = \varphi \Big (x_2 {\partial\over \partial x_1} + x_3 {\partial\over \partial x_2} +\cdots + x_n
{\partial\over \partial x_{n-1}}\Big )
$$
Sachant que :
$$
	\Big ( x_2 {\partial\over \partial x_1} + \cdots + x_n {\partial\over \partial x_{n-1}}\Big )f = \Big
({d\over dt}\Big )_0 \, f(exp(tX)x)
$$
ce r\'esultat signifie que $G$ est caract\'eris\'e par ses polyn\^omes invariants.

\vskip 5mm
\noindent
\textsc{Remarque} : - Supposons $n=3$, et $L = a_1 {\partial\over \partial x_1} + a_2
{\partial\over
\partial x_2} + a_3 {\partial\over \partial x_3}$. Si $L P_2 = 0$, on a : $a_3 = 0$. Si $LP_1 = 0$, on a
:
$$
	a_1 x_3 - a_2 x_2 = 0
$$
d'où $a_1 = \varphi x_2,\ a_2 = \varphi x_3$ et $L = \varphi(x_2 {\partial\over \partial
x_1} + x_3 {\partial\over \partial x_2})$.

\bigskip
- La r\'eduction du cas de la dimension $n$ \`a la dimension $(n-1)$, utilise les propri\'et\'es : 

\medskip
	\quad . La variable $x_1$ n'apparaît que dans le polyn\^ome $P_1$ et ${\partial P_1\over
\partial x_1} = x^{n-1}_n$.

\medskip
	\quad . Les polyn\^omes $P_2,\ldots ,P_n$ jouent le r\^ole de $P_1,\ldots ,P_n$ dans le cas de la
dimension $(n-1)$.

\bigskip
\item Dans le cas $n=2$, il y a un seul polyn\^ome $P(x_1,x_2) = x_2$ et le r\'esultat pr\'ec\'edent n'est
pas valable. On a vu par contre dans 3.2.1. ci-dessus, que l'ensemble des distributions invariantes est
caract\'eristique.

\vskip 4mm
\noindent
\textsc{Conjecture} : Lorsque $G$ est un sous-groupe ferm\'e connexe du sous-groupe unipotent
standard de $GL(n,\Bbb{R})$, l'espace des distributions $G$-invariantes sur $\Bbb{R}^n$ caract\'erise
l'alg\`ebre de Lie des champs de vecteurs associ\'es \`a l'op\'eration de $G$ dans $\Bbb{R}^n$.
\end{enumerate}


\vfill\eject
\renewcommand{\refname}{Bibliographie} 

\end {document}